\documentclass[12pt]{article}
\usepackage{amsxtra,amssymb,amsthm,amsmath,latexsym}

\theoremstyle{plain}
\newtheorem{theorem}{Theorem}

\def\oH{{\overset{\circ}{H}}}
\def\oH1{{\overset{\circ}{H}\kern-.02in{}^1}}

\def\bee{\begin{equation*}}
\def\eee{\end{equation*}}
\def\be{\begin{equation}}
\def\ee{\end{equation}}

\begin{document}
\title{Inverse problems for parabolic equations 2 }

\author{A.G. Ramm\\
 Mathematics Department, Kansas State University, \\
 Manhattan, KS 66506-2602, USA\\
ramm@math.ksu.edu,\\ fax 785-532-0546, tel. 785-532-0580}

\date{}
\maketitle\thispagestyle{empty}

\begin{abstract}
\footnote{MSC:  35K20, 35R30;\, PACS 02.30.Jr}
\footnote{Key words:  parabolic equations, inverse problems, inverse 
source problems}

Let $u_t-u_{xx}=h(t)$ in $0\leq x \leq \pi,\,\,t\geq 0.$ 
Assume that $u(0,t)=v(t)$,  $u(\pi,t)=0$, and $u(x,0)=g(t)$.
The problem is: {\it what extra data determine the three unknown
functions $\{h, v, g\}$ uniquely?}.
This question is answered and an analytical method for recovery 
of the 
above three functions is proposed.

\end{abstract}

\section{Introduction}\label{S:1}

Consider the problem
\be\label{e1}u_t-u_{xx}=h(t)
\quad (x,t)\in 
[0,\pi]\times[0, \,\infty), \ee
\be\label{e2}u(0,t)=v(t), \quad u(\pi, t)=0, \quad u(x,0)=g(x), 
\ee
where the three functions $\{h,\,v,\, g\}$ are not known.

The Inverse Problem (IP) we are interested in is the following one:

{\it What extra data determine the triple   $\{h,\,v,\, g\}$
uniquely?}

There is an extensive literature on inverse problems
for the heat equation (see [1], [2] and references therein),
but the above IP has not been studied, as far as the author knows.
In [3] the author studied an inverse source problem for multidimensional
heat equation in which the source was assumed to be a finite sum of point 
sources, and the inverse problem was to find the location and the 
intensity (strength) of these point sources from experimental data. 
In [4] an inverse problem related to continuation of the solution to heat 
equation is studied.

Let $||f||:=||f||_{L^2(0,\pi)},\,\, u_m:=(u,f_m)=\int_0^{\pi}uf_m dx$,
\be\label{e3}
f_m^{\prime\prime}+m^2f_m=0,\quad 0\leq x \leq \pi,\quad 
f_m(0)=f_m(\pi)=0,\quad 
||f_m||=1,
\quad m=1,2\dots,
\ee
where $f_m=\sqrt{\frac 2 {\pi}}\sin (mx)$. Let $y\in (0,\pi)$ be a point 
such that 
\be\label{e4}
f_m(y)\neq 0\quad \forall m=1,2,\dots
\ee
Our result is:
\begin{theorem}\label{T:1}
The three functions $\{u_1(t),\, u_3(t),\,   
u(y,t)\},$  known for all $t\geq 0$, determine the 
triple $\{h,\,v,\,g\}$  uniquely.
\end{theorem}

We will outline a method for finding 
$h,\,v,\,$ and $g$ and discuss the ill-posedness of the IP.

In Section 2 proofs are given.

\section{Proofs}\label{S:2}

\begin{proof}[Proof of Theorem 1.]

Let us look for the solution to problem (1)-(2) of the form
\be\label{e5}
u(x,t)=\sum_{m=1}^\infty u_m(t)f_m(x),
\ee
where the functions $u_m$ are to be found.
Multiplying equation (1) by $f_m(x)$ and integrating 
over the interval $[0,\pi]$ and then by parts, one gets
\be\label{e6}
\dot {u}_m+m^2u_m=v(t)f'_m(0)+c_m h(t),\quad u_m(0)=g_m,\quad 
c_m:=(1,f_m)=\sqrt{\frac 2{\pi}} \frac {1-\cos(m\pi)}{m},
\ee
where $m=1,2, \dots$.
Thus,
\be\label{e7}
u_m(t)=g_me^{-m^2t}+\int_0^t e^{-m^2(t-s)}[v(s)f'_m(0)+c_mh(s)]ds.
\ee
If the data 
\be\label{e8}
\{u_1(t),\,\, u_3(t), \,\, u(y,t)\}
\ee
are known, then one gets
\be\label{e9}
u_1(t)=g_1e^{-t}+\int_0^t e^{-(t-s)}[v(s)f'_1(0)+c_1h(s)]ds,
\ee
and
\be\label{e10}
u_3(t)=g_3e^{-9t}+\int_0^t e^{-9(t-s)}[v(s)f'_3(0)+c_3h(s)]ds.
\ee
Take $t=0$ in (9) and (10) and get $g_1=u_1(0)$ and $g_3=u_3(0)$.

{\it Thus, $g_1$ and $g_3$ are determined uniquely by the data.} 

Define
$u_1(t)-g_1e^{-t}:=F_1(t),\,\, u_3(t)-g_3e^{-9t}:=F_3(t)$,
and rewrite (9) and (10) as
\be\label{e11}
F_1(t)=\int_0^t e^{-(t-s)}[v(s)f'_1(0)+c_1h(s)]ds,
\ee
and
\be\label{e12}
F_3(t)=\int_0^t e^{-9(t-s)}[v(s)f'_3(0)+c_3h(s)]ds.
\ee
Differentiate (11) and (12) and get 
\be\label{e13}
v(t)f'_1(0)+c_1h(t)=e^{-t}\frac d {dt}[e^{t}F_1(t)]
\ee
\be\label{e14}
v(t)f'_3(0)+c_3h(t)=e^{-9t}\frac d {dt}[e^{9t}F_3(t)]
\ee
This is a linear system for finding $v$ and $h$. The determinant of 
this system is
\be\label{e15}
\left|
\begin{array}{cc}
f_1'(0)\,\, c_1\\
f_3'(0)\,\, c_3\\
\end{array}
\right|=-\frac {32}{3\pi}\neq 0,
\ee
so {\it $v$ and $h$ are uniquely, explicitly and analytically determined 
by the data.} 

If $v(t)$ and $h(t)$ are found, then one has
\be\label{e16} 
u(y,t)=\sum_{m=1}^\infty e^{-m^2t}g_mf_m(y)+w(y,t),
\ee
where $w(y,t)$ is known:
\be\label{e17}
w(y,t)=\sum_{m=1}^\infty f_m(y)\int_0^t 
e^{-m^2(t-s)}[v(s)f_m'(0)+c_mh(s)]ds.
\ee
Denote $q(y,t):=u(y,t)-w(y,t)$. Then $q(y,t)$ is known and
\be\label{e18}
\sum_{m=1}^\infty e^{-m^2t}g_m f_m(y)=q(y,t).
\ee
This relation allows one to determine the numbers 
$g_m f_m(y)$ uniquely for all $m=1,2, \dots$, by the 
formulas:
\be\label{e19}
g_1 f_1(y)=\lim_{t\to \infty}e^tq(y,t), \quad g_2 f_2(y)=\lim_{t\to 
\infty}e^{4t}[q(y,t)-e^{-t}g_1 f_1(y)],
\ee
and so on. Thus, consequitively one finds all the numbers $b_m:=g_m 
f_m(y)$.

If the numbers $b_m$  are found for all $m=1,2, \dots$,  
then the numbers $g_m$ are uniquely determined by the formulas:
\be\label{e20}
g_m=\frac {b_m}{ f_m(y)}.
\ee
Formulas (20) make sense because of the assumption (4).
If all the coefficients $g_m$ are found, then the function $g$
is calculated by the formula:
\be\label{e21}
g(x)=\sum_{m=1}^\infty g_m f_m(x).
\ee
Thus, the triple $\{h,\,v,\,g\}$ is uniquely and analytically found from 
the data \newline
 $\{u_1(t),\,u_3(t),\, u(y,t)\}$, known for all $t>0$. 
Theorem 1 is proved. 
\end{proof}

In the proof of Theorem 1 we assume that the data are exact. The inverse 
problem under discussion is ill-posed: small perturbations of the data
may threw the data out of the set of admissible data. For example, the 
solution $u(x,t)$ is infinitely differentiable (even analytic) with 
respect to $t$ in the region $t>0$, so $u(y,t)$ cannot be an arbitrary 
function. Also, calculation by formulas (19) is an ill-posed problem: 
small errors in calculation of  $g_m f_m(y)$ lead to large errors in
calculation of  $g_{m+1} f_{m+1}(y)$ because of the exponential factor
$e^{-(m+1)^2t}$. A detailed study of a similar problem, arising in
the singularity expansion method (SEM), developed in scattering theory,
is presented in [5], pp.365-393.
Formula (20) also leads to ill-posedness, because the 
denominator in this formula is small for large $m$. Therefore the IP
is severely ill-posed.


\begin{thebibliography}{1000} 


\bibitem{R1} 
Ramm, A.~G.,
{\bf Inverse Problems}, Springer, New York, 2005.

\bibitem{R2}
Ramm, A.~G., {\it An inverse problem for the heat equation},
J.~Math.~Anal.~Appl., 264, N2, (2001), 691-697.

\bibitem{R3}
Ramm, A.~G., {\it Inverse problems for parabolic equations},
Australian Jour. ~Math.~Anal.~Appl., (2005) (to appear).

\bibitem{R4}
Ramm, A.~G., {\it An inverse problem for the heat equation},
 Proc.Roy.Soc. Edinburgh, 123, N6, (1993), 973-976.

\bibitem{R5}
Ramm, A.~G., {\bf Scattering by Obstacles}, D.Reidel, Dordrecht, 1986.

\end{thebibliography}
\end{document}